\documentclass[a4,11pt]{article}
\usepackage{amsmath, amssymb}
\usepackage{graphicx}   
\usepackage{epsfig}     
\usepackage{subfigure}
\newtheorem{proposition}{Proposition}[section]
\newtheorem{lemma}{Lemma}[section]

\newtheorem{remark}{Remark}[section]

\newcommand{\x}{{x}}
\newcommand{\lap}{\bigtriangleup}
\newcommand{\C}{\mathcal C}
\newcommand{\p}{\partial}
\newcommand{\e}{\varepsilon}
\renewcommand{\S} {\mathcal S}
\newcommand{\proof} [1]
   { {\bf Proof.} #1 \hfill\rule{0.5em}{1.2ex} \par\medskip}

\newcommand{\RE}{\mathbb{R}}
\renewcommand{\H}{\mathcal{H}}

\title{Boundary Integral Equations for the Laplace--Beltrami Operator}
\author{S.~Gemmrich$^{1,\natural}$, N.~Nigam$^{1,\sharp}$, O.~Steinbach$^2$}
\date{$^1$Department of Mathematics and Statistics, McGill University, \\
805 Sherbrooke, Montreal H3A 2K6, Canada \\[1mm]
$\sharp$ {\small \tt nigam@math.mcgill.ca}\quad
$\natural$ {\small \tt gemmrich@math.mcgill.ca} \\[2mm]
$^2$Institute of Computational Mathematics, TU Graz, \\
Steyrergasse 30, A 8010 Graz, Austria \\[1mm]
{\small \tt o.steinbach@tugraz.at}}

\begin{document}

\maketitle

\section{Introduction and motivation}
We present a boundary integral method, and an accompanying boundary element discretization, for solving boundary-value problems for the Laplace-Beltrami operator on the surface of the unit sphere $\S$ in $\mathbb{R}^3$.  We consider a closed curve ${\cal C}$ on ${\cal S}$ which divides
${\cal S}$ into two parts ${\cal S}_1$ and ${\cal S}_2$. In particular,
${\cal C} = \partial {\cal S}_1$ is the boundary curve of ${\cal S}_1$. We are interested in solving a  boundary value problem for the Laplace-Beltrami operator in $\S_2$, with boundary data prescribed on $\C$.

We shall begin by describing a physical problem of interest. Then, we derive an integral representation formula for solutions of the Laplace-Beltrami operator on the sphere, and introduce the single and double layer potentials. We investigate their jump properties, and use these to derive an integral equation for the solution of a Dirichlet problem. A variational strategy is presented, along with  some numerical experiments validating our ideas.
To the best of our knowledge, these integral equations and their discretizations have not been studied before. We believe they present an elegant and natural solution strategy for boundary value problems on the sphere.

 This work is motivated in part by recent investigations into the motion of point vortices on spheres, specifically in bounded regions with walls on the sphere. Kidambi and Newton \cite{Kidambi} considered such a problem, assuming the bounded sub-surface of the sphere lent itself to the method of images. Crowdy, in a series of papers \cite{Crowdy2003,Crowdy2004, Crowdy2006}, has also investigated the motion of vortices on spheres. In \cite{Crowdy2006}, he  uses conformal mapping onto the complex plane to study the motion of a vortex on a sphere with walls. We shall study a closely related model problem, for which the methods of \cite{Crowdy2006,Kidambi} would be applicable. However, we propose  an integral-equation method instead which is valid for {\it any} bounded sub-region $\S_2$, provided the curve $\C$ is sufficiently smooth. This technique will be valid even where the method of images is not, and which does not involve explicit knowledge of conformal mappings between the stereographically-projected subregion of interest, and the upper half of the complex plane.

\subsection{Point vortex motion on a sphere with walls}
The underlying physical phenomenon  considered in \cite{Crowdy2006} is the motion of  a point vortex in an incompressible fluid on the surface of the unit sphere, $\mathcal S$. There is a bounded solid region, denoted $\mathcal S_1 \subseteq \mathcal S$, with a simply connected boundary, $\C$.  No fluid can penetrate into $\mathcal S_1$.  Let $\mathcal S_2$ be the surface of the sphere excluding $\mathcal S_1 \cup \C$, see Figure (1).

\begin{figure}[htp]
\centering

\includegraphics[width=0.45\textwidth]{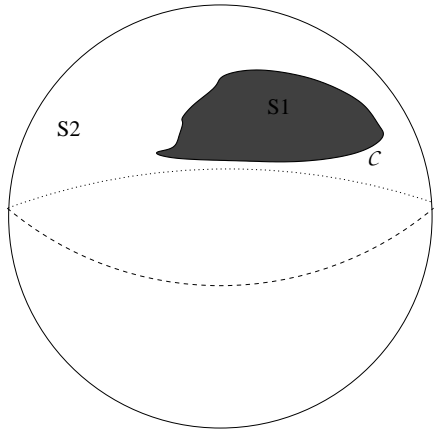} 
\label{fig1}
\caption{The unit sphere, $\S$, with an impenetrable island $\S_1$ on the surface. $\C$ is the boundary of the island.}
\end{figure}
 
A point on the sphere $\S$ will be described in terms of the spherical angles, \[
x(\varphi,\theta) \, = \, \left( \begin{array}{c}
\cos \varphi \sin \theta \\
\sin \varphi \sin \theta \\
\cos \theta
\end{array} \right) \in{\cal S}, \quad
\varphi \in [0,2\pi), \theta \in [0,\pi].
\]
 We consider a point vortex of strength $\kappa$ located at a point $x_0 \in \S_2$. The flow motion is assumed irrotational, except for the point vorticity associated with the vortex. This assumption needs some justification, which we will discuss below.
 The incompressible nature of the fluid allows us to prescribe a stream function, $\Psi({x_0}, {x})$, for the fluid velocity. That is, the velocity field satisfies
 $$ u = \nabla \Psi \times \vec{e}_r.$$ Here $\vec{e}_r$ is the unit radial vector to the surface. The vorticity is then defined as 
 $$\vec{\omega} = \omega \vec{e}_r:=\nabla \times u.$$
If the fluid motion is irrotational except at $\x_0$, then  $\omega=0$ except at that point. We insist that the boundary $\C$ be a streamline of the motion. Without loss of generality, we can set the streamline constant to zero. The function $\Psi$ is really the Green's function for the Laplace-Beltrami operator on the subsurface $\S_2$ of the sphere:
\begin{subequations}
\begin{align}-\lap_\S\Psi(\x_0,\x) &= \kappa\delta(|\x -\x_0|), & \forall \x \in \S_2, \\
\Psi(\x_0,\x) &=0, &\forall \x\in \C.\end{align}\end{subequations}
 We recall that, in spherical coordinates, $\lap_\S$ is defined as
 \begin{equation}\label{LaplaceBeltrami}
\Delta_S u(x) \, = \, \left[
\frac{1}{\sin^2 \theta} \frac{\partial^2}{\partial \varphi^2} +
\frac{1}{\sin \theta} \frac{\partial}{\partial \theta}
\left(
\sin \theta \frac{\partial}{ \partial \theta}
\right)
\right] u(x(\varphi,\theta)).
\end{equation}

 The assumption that the fluid motion is irrotational can be justified by noting that since the fluid is incompressible, it can perfectly slip at $\Sigma$. This allows us to prescribe a circulation at $\C$, so that the Gauss constraint for the vorticity,  
$$ \int_{S_2} \omega \, d\sigma =0, $$ is satisfied.

Analogously to finding Green's functions in the plane, we can find $\Psi$ in terms of the fundamental singularity $U$ of the Laplace-Beltrami operator on the entire surface of the sphere, and a smooth function $v_{x_0}(\x)$. That is, 
\begin{equation} \Psi(\x,\x_0) = U(\x,\x_0) + v_{x_0}(\x)\end{equation}
where  $v_{x_0}$ solves
\begin{equation} \lap_\S v_{x_0}(\x) = \frac{\kappa}{4\pi}, \qquad \forall \x \in \Omega, \qquad v_{x_0}(\x) = -  U(\x_0,\x),\qquad  \x\in \C.\end{equation} 
We can interpret the fundamental singularity, $U$, as the stream function for a point vortex of strength $\kappa$ on the sphere without boundaries. We denote the fundamental singularity, with $\kappa=1$, as $U^*$ henceforth. Note that from \cite{Crowdy2006}, we get that 
$$ U(\x_a,\x)  = -\kappa \log\left\vert\frac{ (z-z_a)(\overline{z-z_a)}}{(1+|z|^2) (1+|z_a|^2)} \right\vert$$ where we've stated this in terms  of the mapped points $z, z_a$ in the complex plane: 
$$ z = \cot(\theta/2) e^{i\phi}, \qquad \x = (\theta, \phi)$$
The fundamental singularity $U$ satisfies the partial differential equation:
\begin{subequations}
\begin{equation} -\lap_\S U = \kappa\left(\delta(|\x-\x_0|) -\frac{1}{4\pi}\right),\label{fund}\end{equation}
and the Gauss condition for the vorticity, $\omega = -\lap_\S U$: 
 \begin{equation} \int_S \omega \,ds =0.\end{equation}\end{subequations}  

The implication of (\ref{fund}) is that there is a "sea" of uniform vorticity, $\frac{1}{4\pi}$, in which the point vortex at $\x_0$ must be embedded if moving on the whole sphere. We cannot find a distribution $\tilde{U}$ which satisfies $ -\lap_S \tilde{U} =\delta{|x-x_0|}$ on the sphere, and which simultaneously satisfies the Gauss constraint. In order to satisfy the constraint and simultaneously have an irrotational flow, we must either counterbalance the point vortex at $x_0$ by another vortex on the sphere, or have the entire fluid moving with a  uniform background vorticity.
 This feature of the fundamental singularity will appear again in the next section, and will require us to impose a side constraint on the solution density, when employing integral equations. This is reminiscent of similar problems arising in the solution of potential problems in unbounded regions of the plane.

At this juncture, we could use the Green's function $\Psi$ to study the motion of the point vortex, which is governed in the stereographic coordinates $z$ by 
$$ \frac{\p z_0}{\p t} = \frac{-i}{2} (1+|z_0|^2)\frac{\p v_{z_0}}{\p z}\vert_{z=z_0}.$$  
The solutions of this evolution equation are described in terms of level sets of the smooth part, $v_{x_0}(x_0)=constant.$  Such an investigation is performed in the papers by Kidambi and Newton \cite{Kidambi}, and Crowdy \cite{Crowdy2006}. Instead, we shall study a closely related mathematical model problem.

Consider the Dirichlet boundary value problem in $\S_2$ for the
Laplace--Beltrami operator: \\
{\it Find a smooth $u$ such that for given Dirichlet data $g$}
\begin{subequations}
\label{DBVP}
\begin{align}
\Delta_{\cal S} u(x) \,& = \, 0 \quad \mbox{for} \; x \in {\cal S_2},\\
u(x) = g(x) \quad \mbox{for} \; x \in {\cal C}
\end{align}\end{subequations}
 We wish to solve (\ref{DBVP}) by reformulating the boundary value problem as an integral equation. As usual, the process of reformulation is not unique; we shall be employing a layer ansatz, and solving an integral equation of the first kind for the unknown density. We note that we could equivalently have chosen to study the Neumann or Robin problem for the system.
We could  use the Green's function for $\S_2$, $\Psi$,  to solve this Dirichlet problem.
The methods suggested in \cite{Kidambi} and \cite{Crowdy2006} would also be applicable for our model problem, with some caveats: the method of Kidambi and Newton relies on the ability to use the method of images, while Crowdy's work requires knowledge of a conformal map from the stereographically-projected $\S_2$ into the upper half plane or the unit circle. Instead, we propose an integral equation method which is valid for {\it any} bounded sub-region $\S_2$, provided the curve $\C$ is sufficiently smooth, and without conformally mapping to the plane. Additionally, integral equations allow us to solve problems with lower regularity properties, a feature we shall explore in upcoming work.

 If $\S_1$ were to degenerate, ie, if the interior of $\C$ had zero area, we would need to add extra conditions to satisfy the Gauss constraint. Mathematically, we would be dealing with the screen problem, and anticipate singular behaviour on the corners of the screen. On  the surface of the entire sphere  without walls, we must either embed the point vortex in a fluid of uniform vorticity (hence no longer irrotational), or counter-balance it by another point vortex.

\section{An integral representation formula on the sphere}
\setcounter{equation}{0}
We begin by reminding the reader of some vectorial identities on the sphere. Let $\vec{e}_r, \vec{e}_\theta, \vec{e}_{\varphi} $ be the usual unit vectors in spherical coordinates. 
Recall that we can define the surface gradient of a scalar $f$ on $\S$ as
\[
\nabla_{\cal S} f(x) \,  = \, \frac{1}{\sin \theta}\frac{\partial f}{\partial \varphi}
\,\vec{e}_{\varphi} + \frac{\partial f}{\partial \theta}\, \vec{e}_{\theta}.
\]
In the same way we introduce the surface divergence for a vector-valued function $\vec{V}$ on the sphere as
\[
\mbox{div}_{\cal S} \vec{V}(x) \, = \, 
\frac{1}{\sin \theta} \left(
\frac{\partial}{\partial \varphi}
V_{\varphi}(\varphi,\theta) +
\frac{\partial}{\partial \theta} ( \sin \theta \, V_{\theta}(\varphi,\theta) )
\right).
\]
We easily see the identity:
\[
\Delta_{\cal S}u(x) \, = \, \mbox{div}_{\cal S} \nabla_{\cal S} u(x) .
\]
We introduce the vectorial surface rotation for a scalar field $f$ on the sphere:
\[
\underline{\mbox{curl}}_{\cal S} f(x) \, = -\, \frac{\partial f}{\partial \theta}
\,\vec{e}_{\varphi} + \frac{1}{\sin \theta}\frac{\partial f}{\partial \varphi}\, \vec{e}_{\theta}
\]
and the (scalar) surface rotation of a vector field $\vec{V}$ as 
\[
\mbox{curl}_{\cal S} \vec{V}(x) \, = \,
\frac{1}{\sin \theta} \left(
- \frac{\partial}{\partial \varphi} V_{\theta}(\varphi,\theta) +
\frac{\partial}{\partial \theta} (\sin \theta \, V_{\varphi}(\varphi,\theta))  
\right).
\]
We then obtain another vectorial identity for the Laplace-Beltrami operator:
\[
\Delta_{\cal S} u(x) \, = \, 
- \mbox{curl}_{\cal S} \underline{\mbox{curl}}_{\cal S} u(x) 
\quad \mbox{for} \; x \in {\cal S}.
\]

We shall be using a variational setting for most of this paper; to this end, we introduce
the inner product
\[
\langle u , v \rangle_{L_2({\cal S})} \, = \, \int\limits_{\cal S}
u(x) v(x) d\sigma_x \, = \,
\int\limits_0^{2\pi} \int\limits_0^\pi
u(x(\varphi,\theta)) v(x(\varphi,\theta)) \sin \theta \, d\theta d\varphi.
\]
We shall now derive the Green's identiy. We find, by  integration by parts,
\begin{eqnarray*}
\langle - \Delta_{\cal S} u , v \rangle_{L_2({\cal S})} \, = \, 
a_{\cal S}(u,v) \, = \, a_{\cal S}(v,u) \, = \,
\langle u , - \Delta_{\cal S} v \rangle_{L_2({\cal S})}
\end{eqnarray*}
where we have introduced the symmetric bilinear form
\begin{align*}
a_{\cal S}(u,v) &:= 
\int\limits_0^{2\pi} \int\limits_0^\pi \left[ 
\frac{1}{\sin \theta}  
\frac{\partial}{\partial \varphi}
u(\varphi,\theta) \frac{\partial}{\partial \varphi}
v(\varphi,\theta) + 
\sin \theta \frac{\partial}{\partial \theta} v(\varphi,\theta) 
\frac{\partial}{\partial \theta}u(\varphi,\theta) \right]
d\theta d\varphi \\
& =  \int\limits_{\cal S} \nabla_S u(x) \cdot \nabla_S v(x) \, d\sigma_x.
\end{align*}

Stoke's theorem for the  positively oriented curve $\C$ and region $\S_2$ may be written  as
\[
\int\limits_{{\cal S}_2} \mbox{curl}_{\cal S} \vec{V}(x) d\sigma_x \, = \,
\int\limits_{\cal C} \vec{V}(x) \cdot \vec{t}(x) \, ds_x .
\] Here, $\vec{t}$ is the unit tangent vector to $\C$. We note that a similar identity holds for the region $\S_1$, with care taken with the orientation of the tangent.
Now, setting $\vec{V} = v(x) \vec{W}(x)$ and applying the product rule we get
\[
\int\limits_{{\cal S}_2} \underline{\mbox{curl}}_{\cal S} v(x) \cdot
\vec{W}(x) \, d\sigma_x \, = \, 
-\int\limits_{\cal C} v(x) [\vec{W}(x) \cdot \vec{t}(x)] ds_x +
\int\limits_{{\cal S}_2} v(x) \mbox{curl}_{\cal S} \vec{W}(x) d\sigma_x .
\]
With $\vec{W}(x) = \underline{\mbox{curl}}_{\cal S} u(x)$ we finally obtain
Green's first formula for the Laplace--Beltrami operator,
\begin{align}
 -\int\limits_{{\cal S}_2} \underline{\mbox{curl}}_{\cal S} v(x) \cdot
 \underline{\mbox{curl}}_{\cal S} u(x) \, d\sigma_x 
=& \int\limits_{\cal C} v(x) [\underline{\mbox{curl}}_S u(x) \cdot 
\underline{t}(x)] ds_x \notag \\
&+\int\limits_{{\cal S}_2} v(x) \Delta_{\cal S}u(x) d\sigma_x.\label{Green 1}
\end{align}
Note that the left hand side of \eqref{Green 1} coincides with the bilinear form 
$a_{\cal S}(u,v)$, with the role of $\S$ being played by $\S_2$.

\subsection{Fundamental solution and a representation formula}

\begin{proposition}{\rm \cite{Firey:1967,Martinez:2005}}
The fundamental solution of the Laplace--Beltrami operator $\Delta_{\cal S}$
as defined in {\rm (\ref{LaplaceBeltrami})} is given by
\begin{eqnarray}
U^*(x,x_0) & = & -
\frac{1}{4\pi} \log |1-(x_,x_0)| \\
&& \hspace*{-3cm} =  - \frac{1}{4\pi} \log \left[
1 - \cos (\varphi-\varphi_0)\sin \theta \sin \theta_0 
- \cos \theta \cos \theta_0 \right].
\end{eqnarray}
In particular,
\begin{equation}
\Delta_{\cal S} U^*(x),x_0) \, = \,
\frac{1}{4\pi} 
\end{equation}
for $x=x(\varphi,\theta),\,x_0= x(\varphi_0,\theta_0) \in {\cal S}$ with
$(\varphi,\theta) \neq (\varphi_0,\theta_0)$.
\end{proposition}

\begin{remark}
For $x, x_0 \in {\cal S}$ we have
\begin{eqnarray*}
|x-x_0|^2 & = &
|x|^2 - 2 (x,x_0)+
|x_0|^2 =
2 [1-(x,x_0)] .
\end{eqnarray*}
Hence we obtain
\[
- \frac{1}{2\pi} \log |x-x_0| \, = \,
- \frac{1}{4\pi} \left[ \log[1-(x,x_0)]
+ \log 2 \right].
\]
In particular, the fundamental solution of the three--dimensional
Laplace--Beltrami operator corresponds to the fundamental solution
of the two-\-dimen\-sional Laplace operator.
\end{remark}

The first Green's identity can be used to derive a representation formula 
for smooth functions defined on the sphere.

\begin{proposition}\label{prop:repr}
Every sufficiently smooth function $u$ on $\S_2$ ($u\in{\cal C}^2({\cal S}_2)\cap{\cal C}^1(\bar{\cal S}_2) $) satisfies the following representation formula
\begin{align}\label{REP}
\displaystyle
  \frac{1}{4\pi}\int\limits_{{\cal S}_2} u(x) d\sigma_x  &-  \int\limits_{{\cal S}_2}
U^*(x,x_0)\Delta_{\cal S}u(x) d\sigma_x  \nonumber \\ 
-  \displaystyle\int\limits_{{\cal C}} U^*(x,x_0) \, \underline{\mbox{curl}}_{\cal S} u(x) \cdot \vec{t}(x) ds_x
 &+ \displaystyle\int\limits_{{\cal C}}  u(x) \,\underline{\mbox{curl}}_{{\cal S}} U^*(x,x_0) \cdot \vec{t}(x) ds_x \nonumber \\
&= \left\{ \begin{array}{ll} u(x_0) & \mbox{ if $x_0\in{\cal S}_1$},\\
				0 & \mbox{ if  $x_0\in{\cal S}\setminus\bar{\cal S}_1$}. 
	\end{array} \right.
\end{align}\end{proposition}
\proof {We obtain Green's second formula by interchanging the roles of $u$ and $v$ in \eqref{Green 1}, adding the two identities and using the symmetry of the left hand side.
\begin{align}
&\displaystyle \int\limits_{{\cal S}_2} u(x)\Delta_{{\cal S}}v(x) - v(x)\Delta_{\cal S}u(x)\, d\sigma_x\notag\\
=&\int\limits_{\cal C} [v(x)\,\underline{\mbox{curl}}_{\cal S} u(x) - u(x)\, \underline{\mbox{curl}}_{\cal S}v(x)]\cdot \vec{t}(x) ds_x.
\end{align}
We define the $\varepsilon$-neighbourhood of $x_0$ on $\cal S$, $B_{\varepsilon}(x_0):=\{y\in{\cal S}:\, |y-x_0| \,>\, \varepsilon\}$ and set ${\cal S}_{2,\varepsilon}:={\cal S}\setminus B_{\varepsilon}(x_0).$ The second Green's formula for ${\cal S}_{2,\varepsilon}$ with $v(x)=U^*(x,x_0)$ yields:
\begin{align}
&\displaystyle \frac{1}{4\pi}\int\limits_{{\cal S}_{2,\varepsilon}} u(x) \,d\sigma_x - \int\limits_{{\cal S}_{2,\varepsilon}} U^*(x,x_0)\,\Delta_{\cal S} u(x) \,d\sigma_x\notag\\
=&\int\limits_{\cal C} [ U^*(x,x_0)\,\underline{\mbox{curl}}_{\cal S} u(x) - u(x)\,\underline{\mbox{curl}}_{\cal S} U^*(x,x_0)]\cdot \vec{t}(x) ds_x\notag\\
+&\int\limits_{\partial B_{\varepsilon}(x_0)} [ U^*(x,x_0)\,\underline{\mbox{curl}}_{\cal S} u(x) - u(x)\,\underline{\mbox{curl}}_{\cal S} U^*(x,x_0)]\cdot \vec{t}(x) ds_x\label{Green_S1eps}
\end{align}
First we observe that
\[ |\displaystyle \int\limits_{{\cal S}_{2}} U^*(x,x_0)\,\Delta_{\cal S} u(x) \, d\sigma_x|\, \leq \,
||\Delta_{\cal S}\,u||_{L^{\infty}({\cal S})} \int\limits_{{\cal S}_{2}} |U^*(x,x_0)| \, d\sigma_x\,\leq M,\]
and hence:
\[\displaystyle \lim_{\varepsilon\rightarrow 0} \int\limits_{{\cal S}_{2,\varepsilon}} U^*(x,x_0)\,\Delta_{\cal S} u(x) \,d\sigma_x = \int\limits_{{\cal S}_2} U^*(x,x_0)\,\Delta_{\cal S} u(x) \,d\sigma_x. \]
Furthermore, we can estimate the integral 
\[\left|\int\limits_{\partial B_{\varepsilon}} U^*(x,x_0)\,\underline{\mbox{curl}}_{\cal S} u(x)\cdot\vec{t}(x)\, ds_x\right|\, 
\leq\, ||\underline{\mbox{curl}}_{\cal S}u||_{L^{\infty}} \, \int\limits_{\partial B_{\varepsilon}} |U^*(x,x_0)|\,ds_x,\]
By changing coordinates, without loss of generality $x_0$ can be taken to be the north pole of the sphere, i.e. $x_0=\left(0,0,1\right)^\top$. The curve $\partial B_{\varepsilon}$ is then fully described by the latitude $\theta_{\varepsilon}$, say. According to the cosine law we have $\cos\theta_{\varepsilon}=1-\frac{\varepsilon^2}{2}$ and thus
\begin{align*}
\int\limits_{\partial B_{\varepsilon}}|U^*(x,x_0)|\, ds_x &=\frac{1}{4\pi}\int_{2\pi}^0 \log|1-\cos\theta_{\varepsilon}| \,\sin\theta_{\varepsilon} \, d\varphi\\
&= -\frac{1}{2}\varepsilon \sqrt{(1-\frac{\varepsilon^2}{4})}\,\log\frac{\varepsilon^2}{2}
\longrightarrow 0 \quad \mbox{ (as $\varepsilon\rightarrow 0$)}.
\end{align*}
To analyse the second contribution along $\partial B_{\varepsilon}$, we again assume $x_0$ to be the northpole. We then  compute 
\[\underline{\mbox{curl}}_{\cal S} U^*(x,x_0)=\frac{\sin\theta_{\varepsilon}}{4\pi (1-\cos\theta_{\varepsilon})} \vec{e}_\varphi.  \]
Since the line element on the surface of the sphere is given by
\[\vec{t}(x(\varphi,\theta))\cdot ds_{x(\varphi,\theta)}=  d\theta\, \vec{e}_{\theta} \,+\, \sin\theta d\varphi \,\vec{e}_{\varphi}, \]
we deduce that (note the orientation of $\partial B_{\varepsilon}$):
\begin{align*}
-\int\limits_{\partial B_{\varepsilon}} u(x) [\underline{\mbox{curl}}_{\cal S}U^*(x,x_0)\cdot \vec{t}(x)] ds_x 
&=-\frac{\sin^2\theta_{\varepsilon}}{4\pi(1-\cos\theta_{\varepsilon})}\int_{2\pi}^0 u(\varphi,\theta_{\varepsilon})\, d\varphi, 
\intertext{which in the limit as $\varepsilon\longrightarrow 0$ tends to}
\displaystyle \frac{u(x_0)}{2}\, \lim_{\varepsilon\rightarrow 0}\frac{\varepsilon^2\, (1-\frac{\varepsilon^2}{4})}{\frac{\varepsilon^2}{2}}\,=\, u(x_0).
\end{align*}
Taking the limit as $\epsilon \rightarrow 0$ in \eqref{Green_S1eps} proves the result.
}

At this juncture, we draw the reader's attention to the term $\int_{\S_2} u(x) \,d\sigma_x$ in the representation formula (\ref{REP}). If $u$ satisfied $\lap_\S u=0$ in $\S_2$, the familiar integral representation formula for the Laplacian in 2-D would not involve such a term; indeed, the left hand side of the representation formula consists of line integrals only if $u$ satisfies the side constraint, $\int_{\S_2} u\,d\sigma_x=0$. This is linked to the Gauss constraint.
\section{Layer potentials and boundary integral operators}

Having derived an integral representation formula for solutions of the Laplace-Beltrami problem in $\S_2$ in the previous section, we are now in a position to reformulate the boundary value problem as an integral equation.
\setcounter{equation}{0}

\subsection{Single and double layer potentials}
Following the integral representation derived in Proposition \ref{prop:repr} we
define the following two layer potentials:
\begin{itemize}
\item The {\bf single layer potential} with sufficiently smooth density function $\sigma$:
	\[(\widetilde{V}\sigma)(x):=\int\limits_{{\cal C}} U^*(x,y)\,\sigma(y)\,ds_y\quad\mbox{for $x\not\in{\cal C}$},\]
\item and the {\bf double layer potential} with sufficiently smooth density function $\mu$:
	\[(\widetilde{W}\mu )(x):=\int\limits_{{\cal C}} \mu(y)\,[\underline{\mbox{curl}}_{\cal S} U^*(x,y)\cdot \vec{t}(y)]\, ds_y \quad\mbox{for $x\not\in{\cal C}$}.\]
\end{itemize}
By Proposition  \ref{prop:repr}, every solution to the homogeneous Laplace-Beltrami equation can be written as the sum of a single and a double layer potential modulo a constant. This  is the starting point for the so-called direct boundary integral approach. However, for the purpose of this paper we follow the layer ansatz based on the following observation.

For $x\not\in{\cal C}$, the single layer potential satisfies:
\begin{subequations}\begin{align}
\Delta_{\cal S} (\widetilde{V}\sigma)(x) & =  \Delta_{\cal S}
\int\limits_{\cal C} U^*(x,y)\, \sigma(y)\,ds_{y} 
 =  \int\limits_{\cal C} \Delta_{\cal S} U^*(x,y)\, \sigma(y)\,ds_{y} \notag \\ 
& = \frac{1}{4\pi}\,\int\limits_{\cal C} \sigma(y)\,ds_{y} =  \, 0  \\
\mbox{under the constraint }
&\displaystyle\int\limits_{\cal C} \sigma(y)\,ds_{y} \, = \, 0. \label{SL_constraint}
\end{align}\end{subequations}
Hence, we may find the general solution of the
Dirichlet boundary value problem \eqref{DBVP} as
\begin{align}\label{solution}
u(x) & =  \int\limits_{\cal C} U^*(x,y) \,\sigma(y) \,ds_{y} + p, 
\end{align}
where $p \in \mathbb{R}$ is some Lagrange multiplier related
to the constraint \eqref{SL_constraint}. 
Similarly, the double layer potential satisfies the Laplace-Beltami equation for $x\not\in{\cal C}$
\begin{align}
\Delta_{\cal S} (\widetilde{W}\mu)(x) & = \Delta_{\cal S}
\int\limits_{\cal C} \mu(y)\, [\underline{\mbox{curl}}_{\cal S} U^*(x,y)\cdot\vec{t}(y)]\,ds_{y} \notag\\
&=\int\limits_{\cal C} \mu(y) [\Delta_{\cal S}\,\underline{\mbox{curl}}_{\cal S}\,  U^*(x,y)\cdot\vec{t}(y)]\,ds_{y} \notag \\ 
&=\int\limits_{\cal C} \mu(y) [\underline{\mbox{curl}}_{\cal S} \, \Delta_{\cal S} U^*(x,y)\cdot\vec{t}(y)]\,ds_{y} \notag \\ 
&=0,\notag
\end{align}
without any further constraints on the density $\mu$.
We might thus also try to look for the solution to \eqref{DBVP} in the form of a double layer.

\subsection{Jump relations for $\widetilde{V}$ and $\widetilde{W}$}
In the previous section, we have only defined the layer potentials for $x$ away from the boundary curve. However, in order to align the operators with the given Dirichlet data along $\C$, we need to investigate their behavior in the limit as $x$ approaches $\C$. 
Similarly, if one is interested in solving the Neumann problem in which the tangential component of the vectorial surface rotation is prescribed along $\C$, one has to investigate the limit features of this quantity for the layer potentials. In both cases, there will be certain jump relations across the curve $\C$. For the purpose of this paper however, we will restrict ourselves to the Dirichlet case.
First, consider 
the single layer potential with density $\sigma$ for $\tilde{x} \not \in {\cal C}$:
\begin{align}
(\widetilde{V}\sigma)(\tilde{x}) & =  \int\limits_{\cal C} U^*(\tilde{x},x)\,\sigma(x)\,ds_{x} \notag\\
& = - \frac{1}{4\pi} \int\limits_{\cal C} \log [1-\langle \tilde{x}\,,\,x\rangle]\,\sigma(x)\,ds_{x}
\end{align}
The following lemma describes the limit behavior of the single layer potential.
\begin{lemma}
For $x_0\in{\cal C}$ we have:
\begin{align*} 
(V\sigma)(x_0) & := \lim\limits_{{\cal S}\ni\tilde{x}\to x_0}\,(\widetilde{V}\sigma)(\tilde{x}) = \int\limits_{\cal C} U^*(x_0,y)\,\sigma(y)\,ds_{y}
\end{align*}
as a weakly singular line integral and hence $(\widetilde{V}\sigma)$ is continuous across ${\cal C}$.
\end{lemma}
\proof{Fix an arbitrary $\varepsilon > 0$. Let
$x_0 \in {\cal C}$ be fixed, and  and
$\tilde{x} \in {\cal S}$ satisfy
$|\tilde{x}-x_0| < \varepsilon$. Introduce the notation
$$ \C_{\e,\leq}:=\{y\in \C,|y-x_0|\leq\e\},\qquad  \C_{\e,>}:=\{y\in \C,|y-x_0|>\e\}.$$

Then, if we define
\begin{align*}
I_{\varepsilon}(\tilde{x}) &: = \int\limits_{\cal C} U^*(\tilde{x},y)\,\sigma(y) \,ds_{y} 
 - \int_{\C_{\e,>}} U^*(x_0,y)\,\sigma(y)\,ds_{y}, \end{align*} we can easily show

\begin{align} I_{\varepsilon}& =  \int_{\C_{\e,\leq}}
 \left[ U^*(\tilde{x},y)\,-\, U^*(x_0,y) \right]\,\sigma(y)\,ds_{y} +  \int_{\C_{\e,\leq}}
 U^*(\tilde{x},y)\, \sigma(y)\, ds_{y}.\label{lemm}
\end{align}
The first integral in (\ref{lemm}) vanishes in the limit as $\tilde{x}\longrightarrow x_0$, i.e.
\begin{align*}
\lim\limits_{\tilde{x} \to x_0}
\int_{\C_{\e,>}}
\left[U^*(\tilde{x},y)\,-\,U^*(x_0,y) \right]\,\sigma(y)\,ds_{y} = 0 .
\end{align*}
The second term in (\ref{lemm}) we can bound in terms of the density $\sigma$:
\begin{align*}
 &\left|\int_{\C_{\e,\leq}}
U^*(\tilde{x},y)\,\sigma(y)\,ds_{y}\,\right| \, \leq \,
\| \sigma \|_{L_\infty({\cal C})} \int_{\C_{\e,\leq}}
\left| U^*(\tilde{x},y)\right|\,ds_{y}.
\end{align*}
To finish the proof, note that we can estimate
\begin{align*} \int_{\C_{\e,\leq}}
\left| U^*(\tilde{x},y)\right|\,ds_{y}\,\leq\,
\hspace*{-8mm} \int\limits_{\footnotesize
\begin{array}{c}y \in {\cal C} \\|y-\tilde{x}| \leq 2\varepsilon\end{array}} \hspace*{-6mm}
\left| U^*(\tilde{x},y)\right|\,ds_{y}
\buildrel {\tilde{x}\rightarrow x_0}\over\longrightarrow 
\hspace*{-8mm} \int\limits_{\footnotesize
\begin{array}{c}y \in {\cal C} \\|y-x_0| \leq 2\varepsilon\end{array}} \hspace*{-6mm}
\left| U^*(x_0,y)\right|\,ds_{y}
\buildrel{\varepsilon\rightarrow 0}\over\longrightarrow 0.
\end{align*}
Putting these estimates together, we see that  $\displaystyle \lim_{\varepsilon\rightarrow 0}\lim_{\tilde{x}\rightarrow x} I_{\varepsilon}(\tilde{x})=0$, which proves the assertion.
}

The case of the double layer potential is slightly more involved, since the limit process reveals a hidden delta function. To see this, consider the double layer potential with density $\mu$ for $x\not\in {\cal C}$:
\begin{align}
(\widetilde{W}\mu)(x) & = 
\int\limits_{\cal C} \mu(\tilde{x})\, 
\left[\underline{\mbox{curl}}_{\cal S} 
U^*(x,\tilde{x}) \cdot \vec{t}(\tilde{x})\right]\,
ds_{\tilde{x}} \notag\\ 
 &\hspace*{-14mm}= \frac{1}{4\pi}\int\limits_{\cal C} \mu(\tilde{x})\,\frac{1}{A(x,\tilde{x})}\left[\left(\begin{array}{c}
-\cos (\varphi-\widetilde\varphi)\cos\theta\sin\widetilde\theta+\sin \theta\cos\widetilde\theta\\
-\sin\widetilde\theta\sin (\varphi-\widetilde\varphi)\end{array}\right)
 \cdot \vec{t}(\tilde{x})\right]\, ds_{\tilde{x}}
\end{align}
where $x=x(\varphi,\theta)$ and $\tilde{x}=x(\widetilde\varphi,\widetilde\theta)$ and 
$A(x,\tilde{x})=1-\cos(\varphi-\widetilde\varphi)\sin\theta\sin\widetilde\theta-\cos\theta\cos\widetilde\theta$.

\begin{lemma}
For $x_0\in{\cal C}$ we have:
\begin{align*} 
&(\gamma^{{\cal S}_2}_0\widetilde W\mu)(x_0):=\displaystyle\lim\limits_{{\cal S}_2\ni x\to x_0}\,(\widetilde{W}\mu)(x) \\
&\hspace*{5mm} = (K\mu)(x_0)\,+\,\left(1-\frac{\alpha(x_0)}{2\pi}\right)\,\mu(x_0),
\end{align*}
where $\alpha(x_0)$ represents the interior (with respect to ${\cal S}_2$) angle of ${\cal C}$ at $x_0$. For a smooth curve, $\alpha = \pi$. The operator $(K\mu)$ is given by the following integral expression:
\begin{align}
(K\mu)(x_0)&=\displaystyle\lim_{\varepsilon\rightarrow 0}\, \hspace*{5mm}(K_{\varepsilon}\mu)(x_0)\notag\\
&= \displaystyle\lim_{\varepsilon\rightarrow 0} \hspace*{-5mm}\int\limits_{\footnotesize
\begin{array}{c}\tilde{x} \in {\cal C} \\|\tilde{x}-x_0| < \varepsilon\end{array}} \hspace*{-5mm}
\mu(\tilde{x})\,\left[\underline{\mbox{curl}}_{\cal S}\,U^*(x,\tilde{x})\cdot\vec{t}(\tilde{x})\right]\,ds_{\tilde{x}}.
\end{align}
Hence  the double layer potential satisfies:
\begin{align}
\left[(\gamma_0\widetilde W\mu) \right]_{\cal C}:=(\gamma^{{\cal S}_2}_0\widetilde W\mu)+(\gamma^{{\cal S}_1}_0\widetilde W\mu)=\mu,
\end{align}
where we tacitely assumed the orientation of the tangential vector $\vec{t}$ along ${\cal C}$ to be in accordance with the orientation of ${\cal S}_2$ in the sense of Stoke's theorem.
\end{lemma}

\proof{Given $\varepsilon \, > \,0$, let $x\in {\cal S}_2$ with $ \| x-x_0\|<\varepsilon$.We introduce the notation $$\C_{\e,<}:= \{\tilde{x} \in {\cal C}, |\tilde{x}-x_0| < \varepsilon\},\qquad \C_{\e,\geq}:= \{\tilde{x} \in {\cal C}, |\tilde{x}-x_0| \geq \varepsilon\}$$

 Then,
\begin{align*}
&(W\mu)(x)-(K_{\varepsilon} \mu)(x_0) &&=  \int_{\C_{\e,\geq}}
\mu(\tilde{x})\left[\underline{\mbox{curl}}_{\cal S}\, U^*(x,\tilde{x})\,-\,\underline{\mbox{curl}}_{\cal S}\,U^*(x_0,\tilde{x})\right]\cdot\vec{t}(\tilde{x})\, ds_{\tilde{x}}\\
&&&\int_{\C_{\e,<}}
\mu(\tilde{x})\,\underline{\mbox{curl}}_{\cal S}\,U^*(x,\tilde{x})\cdot\vec{t}(\tilde{x})\,ds_{\tilde{x}}
\end{align*}
and the first integral again vanishes as $x$ approaches $x_0$, i.e.
\begin{align*}
\left| \int_{\C_{\e,\geq}}
\mu(\tilde{x})\left[\underline{\mbox{curl}}_{\cal S}\, U^*(x,\tilde{x})\,-\,\underline{\mbox{curl}}_{\cal S}\,U^*(x_0,\tilde{x})\right]\cdot\vec{t}(\tilde{x})\, ds_{\tilde{x}}\right| \;
\buildrel{x\rightarrow x_0}\over\longrightarrow \; 0.
\end{align*}
The second term can be rewritten as follows:
\begin{align}
&\int_{\C_{\e,<}}
\mu(\tilde{x})\,\underline{\mbox{curl}}_{\cal S}\,U^*(x,\tilde{x})\cdot\vec{t}(\tilde{x})\,ds_{\tilde{x}} \nonumber\\
&=\hspace{8mm}\int_{\C_{\e,<}}
\left[\mu(\tilde{x})-\mu(x_0)\right]\,\underline{\mbox{curl}}_{\cal S}\,U^*(x,\tilde{x})\cdot\vec{t}(\tilde{x})\,ds_{\tilde{x}}\nonumber\\
&\hspace{4mm}+\hspace{8mm}
 \mu(x_0)\int_{\C_{\e,<}}
\underline{\mbox{curl}}_{\cal S}\,U^*(x,\tilde{x})\cdot\vec{t}(\tilde{x})\,ds_{\tilde{x}}. \label{interm}
\end{align}
For the first integral on the right hand side of  (\ref{interm}) we have the estimate
\begin{align*}
& \left|\int_{\C_{\e,<}}
\left[\mu(\tilde{x})-\mu(x_0)\right]\,\underline{\mbox{curl}}_{\cal S}\,U^*(x,\tilde{x})\cdot\vec{t}(\tilde{x})\,ds_{\tilde{x}}\right| \\
&\leq \sup_{\C_{\e,<}} |\mu(\tilde{x})-\mu(x_0)|
\int_{\C_{\e,<}}
\left|\underline{\mbox{curl}}_{\cal S}\,U^*(x,\tilde{x})\cdot\vec{t}(\tilde{x})\right|\,ds_{\tilde{x}}\\
&\leq\hspace*{6mm} M\,\cdot \mbox{length $({\cal C}_{\varepsilon})$}\;\cdot
\sup_{\C_{\e,<}} |\mu(\tilde{x})-\mu(x_0)|
\end{align*}
for some constant $M$, and hence the integral vanishes in the limit as $\varepsilon\longrightarrow 0$.
For the second integral in (\ref{interm}),  we define $\Omega_{\varepsilon}(x):=\left\{ \tilde{x}\in{\cal S}_2:\; |x-\tilde{x}|< \varepsilon \right\}$ to see
\begin{align*}
 \mu(x_0)\int_{\C_{\e,<}}
\underline{\mbox{curl}}_{\cal S}\,U^*(x,\tilde{x})\cdot\vec{t}(\tilde{x})\,ds_{\tilde{x}} 
&= \mu(x_0)\hspace*{-3mm}\int\limits_{\footnotesize
\partial \Omega_{\varepsilon}(x_0)} \hspace*{-3mm}
\underline{\mbox{curl}}_{\cal S}\,U^*(x,\tilde{x})\cdot\vec{t}(\tilde{x})\,ds_{\tilde{x}}\\
&  \hspace*{-10mm}- \mu(x_0)\hspace*{-6mm}\int\limits_{\footnotesize
\begin{array}{c}\tilde{x} \in {\cal C} \\|\tilde{x}-x_0| = \varepsilon\end{array}} \hspace*{-8mm}
\underline{\mbox{curl}}_{\cal S}\,U^*(x,\tilde{x})\cdot\vec{t}(\tilde{x})\,ds_{\tilde{x}} 
\intertext{Using the representation formula with $u\equiv1$ we get}
 =\; \mu(x_0)\left(1-\frac{1}{4\pi}\int\limits_{\Omega_{\varepsilon}}d\sigma_{\tilde{x}}\right) 
& - \hspace*{2mm}\mu(x_0)\hspace*{-6mm}\int\limits_{\footnotesize
\begin{array}{c}\tilde{x} \in {\cal C} \\|\tilde{x}-x_0| = \varepsilon\end{array}} \hspace*{-8mm}
\underline{\mbox{curl}}_{\cal S}\,U^*(x,\tilde{x})\cdot\vec{t}(\tilde{x})\,ds_{\tilde{x}} 
\end{align*}
and without loss of generality we compute the remaining integral with respect to the northpole to find for all $x_0$: 
\begin{align*}
\displaystyle\lim_{\varepsilon\rightarrow 0} \hspace*{-6mm}\int\limits_{\footnotesize
\begin{array}{c}\tilde{x} \in {\cal C} \\|\tilde{x}-x_0| = \varepsilon\end{array}} \hspace*{-8mm}
\underline{\mbox{curl}}_{\cal S}\,U^*(x,\tilde{x})\cdot\vec{t}(\tilde{x})\,ds_{\tilde{x}} =\frac{\alpha(x_0)}{2\pi}.
\end{align*}
Putting the parts together we see that
\begin{align*}
\displaystyle\lim_{\varepsilon\rightarrow 0} \lim_{x\rightarrow x_0} \left((\widetilde{W}\mu)(x)-(K_{\varepsilon}\mu)(x_0)\right)
=\left(1-\frac{\alpha(x_0)}{2\pi}\right)\mu(x_0).
\end{align*}
}

\section{A BIE strategy for solving the Dirichlet problem}
\setcounter{equation}{0}
With the single and double layer potentials defined as in the previous section, we are now in a position to reformulate the Dirichlet problem for the Laplace-Beltrami operator.  For the purposes of this paper, we assume sufficient smoothness of the data $g$ and the curve $\C$ such that all the operators are well-defined; these assumptions can be relaxed, and  the precise regularity and smoothness assumptions necessary are a subject of a forthcoming work. In what follows, however, we assume the curve $\C$ is at least $C^2$, and that the boundary data is smooth. We shall present a numerical example where we allow $g$ to be Lipschitz; 

Recall that we wish to find a smooth function $u$ such that
\begin{equation}-\lap_\S u=0\,\,\mbox{in}\,\,\S_2, \,\,\mbox{and}\,\, u=g\,\,\mbox{on} \,\C. \label{PDEagain} \end{equation}
We seek a solution of this equation in terms of a layer ansatz. That is, we wish to find a density, $\sigma$ or $\mu$, so that either  \begin{equation} u:=\widetilde{V}\sigma  \label{layeransatz}  \,\,\mbox{ or } \,\, u:=\widetilde{W} \mu\end{equation} solves \ref{PDEagain}.

\begin{lemma} If the density $\sigma$ solves the boundary integral equation
\begin{align} V \sigma =g, \,\,\mbox{on}\,\,  \C, \,\mbox{and also} \; \int_{\C} \sigma ds =0,\end{align}
then 
the function $u:= \widetilde{V} \sigma$ solves (\ref{PDEagain}).If the density $\mu$ solves 
\begin{align} (\frac{1}{2} I + K )\mu=g, \qquad \mbox{on}\, \C ,\end{align} then the solution of (\ref{PDEagain}) is given by the double layer potential, $ u:=\widetilde{W}\mu.$
\end{lemma}

The proof is immediate from the previous section. 

\subsection{An indirect integral equation formulation}
For concreteness, we describe in some detail a variational strategy to solve a boundary integral equation, whose solution then can be used to solve (\ref{PDEagain}). We seek a solution $u$ of the Laplace-Beltrami operator with prescribed boundary values on $\C$. The solution is assumed to be of the form
$$ u = \widetilde{V}\sigma + p$$ where the density satisfies the integral equation
\[ V\sigma \,+p\,= g \,\, \mbox{on}\,\,\C\] along with the constraint $$ \int_\C \sigma ds = 0.$$
We can write the weak formulation of this problem in saddle-point form as: {\it Find $\sigma \in \H$ and a multiplier $p\in \mathbb{R}$ such that}
\begin{subequations}
\begin{align}\label{eqn:saddle}
 \langle V\sigma, \chi \rangle + p \langle 1, \chi\rangle &= \langle g,\chi\rangle, \\
 q \langle 1,\sigma\rangle &=0,\end{align}
 \end{subequations}
for any  test function $\chi \in \H$ and any real constant $q$. Under the present assumptions on smoothness, we set $\H= C(\C)$, and $\langle \cdot, \cdot\rangle$ is simply the $L^2-$ inner product along $\C$; we will describe the appropriate Sobolev spaces in which to naturally seek $\sigma$ in subsequent work. The angle brackets will then represent the duality pairings in $L^2$.
 
 The discretization strategy is now standard. Let $\tau_h$ be a partition of $\C$, with sub-interval size $h>0$. We approximate $\H$ by a finite-dimensional space, $S_h$, which is parametrized by the meshsize $h$; as $h\rightarrow 0$, the approximation error $\inf_{v_h\in S_h} \|u-v_h\|_\H \rightarrow 0$ for all $u \in \H$. We then study the discrete Galerkin problem: \\
 {\it Find $\sigma_h \in S_h,\, p_h \in \RE,$ such that for all $(\chi_h, q_h) \in S_h \times \RE$, }
 \begin{subequations}\label{galerkin}
 \begin{align}
 \langle V\sigma_h, \chi_h \rangle + p_h\langle 1, \chi_h\rangle ds &= \langle g,\chi_h\rangle, \\
 q_h \langle 1,\sigma_h\rangle &=0,\end{align}
\end{subequations}
 We shall provide an error analysis of this system in a subsequent paper, based on the correct choices of Sobolev spaces for the densities $\sigma$ and approximation spaces $S_h$; at present, we present numerical experiments to validate the boundary element strategy.

\subsection{Numerical experiments}
In what follows, we choose $\C$ to be the equator of the sphere, which is described by the latitude $\theta = \pi/2$.  We solve the Laplace-Beltrami equation in the southern hemisphere ${\cal S}_2$. To do this we prescribe Dirichlet data $g$ on $\C$ and solve the discrete Galerkin system (\ref{galerkin}). 
The partitions $\tau_h = \cup_{i=1}^N \Omega_i$ of $\C$ are chosen to consist of uniform sub-intervals of size $h$. The approximations are sought in the space of piece-wise constant functions, i.e.
$$ S_h :=\{ \chi \, \vert\, \chi(\varphi)=c_i \;\mbox{for}\;\varphi \in \Omega_i\}.$$
In the specific case of the southern hemisphere, a Green's function for the problem is known and we can write down the solution in closed form as follows: 
\begin{align}
u(\varphi,\theta)&=\frac{1}{4\pi}\int_0^{2\pi} \frac{\cos(\theta)}{-1 + \cos(\varphi-\varphi_0)\sin(\theta)}\, g(\varphi_0)\,d\varphi_0,
\end{align}
for $\theta >\frac{\pi}{2}$  and $\varphi\in [0,2\pi]$. This expression serves as a reference for our computed solution.

\begin{figure}[htp]
\centering
\subfigure[Convergence behavior for different choices of $g$]{\label{fig:conv}
\includegraphics[width=0.45\textwidth]{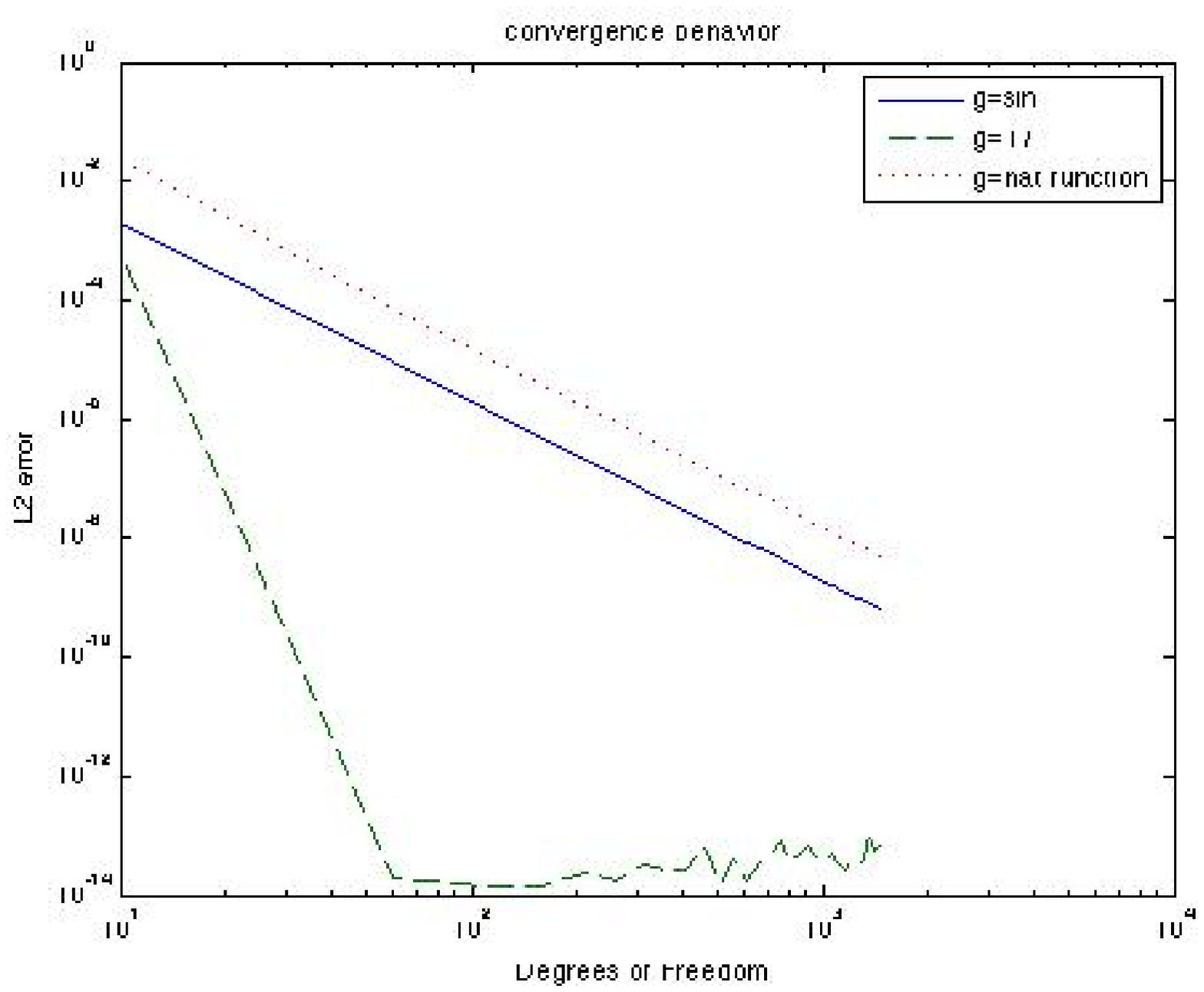} }
\subfigure[View on solution from below the south pole for hat-shaped data, $h(\varphi)$. ]{
\label{fig:bottom}
 \includegraphics[width=0.45\textwidth ]{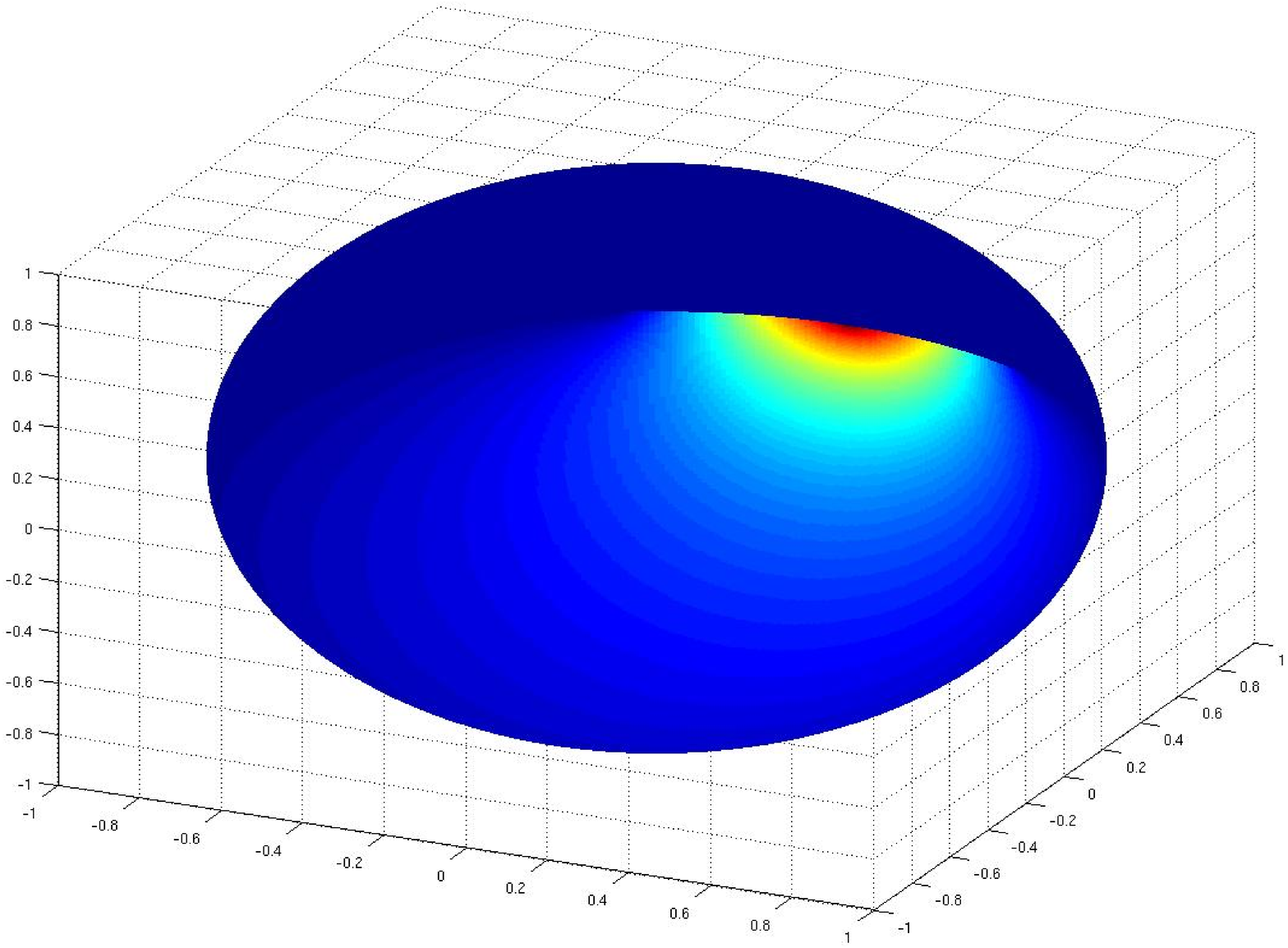}}
\end{figure}

We report the convergence behavior of the method in terms of the $L^2$ error of the computed solution along the latitude $\theta=3/4\,\pi$. In Table 1, we report this $L^2$-error versus the number of unknowns, for two choices of Dirichlet data: $g=\sin(\varphi)$, and $g=h(\varphi)$, where $h$ is a hat-shaped function, $$ h(\varphi):=\begin{cases}
17\left(1-\frac{6}{\pi}|\varphi|\right), & |\varphi|\leq \frac{\pi}{6}\\
0& otherwise. \end{cases} $$
We note that halving the mesh-size reduces the $L^2$ error by a factor of 8 in both cases.
\begin{table}
\caption{$L^2$ error of BEM solution, measured along $\theta=3/4\,\pi$}
\begin{tabular}{c|c|c||c|c}

& \multicolumn{2}{c||}{Convergence for $g=\sin(\varphi)$}& \multicolumn{2}{c}{Convergence for $g=h(\varphi)$ } \\ \hline\hline
DoF&$L^2$ error & Ratio& $L^2$ error & Ratio \\ \hline
20&2.41e-4& -- & 0.023&-- \\ \hline 

40 & 2.99e-5 & 8.06 & 2.57e-4 & 89.49\\ \hline
80& 3.73e-6  & 8.02 & 3.06e-5 & 8.40 \\ \hline
160& 4.66e-7&8.00 &   3.73e-6&  8.20 \\ \hline
320& 5.83e-8&7.99 &   4.61e-7&  8.09 \\ \hline
640&7.28e-9&8.01 &    4.73e-8& 8.05 \\ \hline
\end{tabular}
\end{table}
 Figure \ref{fig:conv} shows the convergence behavior in terms of the above mentioned error versus the number of unkowns. In figure \ref{fig:bottom}, we see the actual solution of the boundary value problem,  corresponding to the piece-wise Dirichlet data $g=h(\varphi)$.

\section{Conclusion}
We have presented a boundary integral formulation, and associated Galerkin discretization strategy, for a boundary value problem for the Laplace-Beltrami operator on the unit sphere in $\RE^3$. Numerical experiments verify the applicability of the idea, and a rigorous error analysis will be presented in future work.

\section{Acknowledgements}
The authors wish to thank the organizers of the Abel Symposium 2006. NN was supported by NSERC and FQRNT.

\end{document}